# A Global root bracketing method with adaptive mesh refinement

*Mohammad Amin Razbani[1]*

*Department of Chemical Engineering, Ferdowsi University of Mashhad, Mashhad, Iran.*

*Email: amin.razbani@gmail.com*

Abstract: *An efficient method for finding all real roots of a univariate function in a given bounded domain is formulated. The proposed method uses adaptive mesh refinement to locate bracketing intervals based on bisection criterion for root finding. Each bracketing interval encloses one root. An adaptive form of error is introduced to enclose roots in a desired tolerance based on how much close the roots are. Detecting roots with even multiplicity, which is regarded out of the realm of bracketing methods, becomes possible with the method proposed in this paper. Also, strategies for finding odd-multiple roots with the least number of function evaluations is proposed. Adaptive mesh refinement lead to considerable reduction in function evaluations in comparison to previous global root bracketing methods. The reliability of the proposed method is being illustrated by several examples.*

## 1    Introduction

A one dimensional root finding problem finds x such that $f(x) = 0$, for a given function $f: R \to R$. Such an x is called a root or zero of the function *f*. Numerical method algorithms which deal with solving this problem can be divided into two basic groups, bracketing methods and open methods. Bracketing methods start with a bounded interval which is guaranteed to bracket a root. The size of initial bracket is reduced step by step until it encloses the root in a desired tolerance. Open methods begin with an initial guess of the root and then improving the guess iteratively. Bracketing methods provide an absolute error estimate on the root's location and always work but converge slowly. In contrast, open methods do not always converge. There is a trade-off between absolute error bound and speed in open methods. Bracketing methods use function evaluations but open methods need both the function and its derivative.

There are three major kinds of bracketing methods. The first one is bisection method. The second one is false position which is a method of finding roots based on linear interpolation. The third one is the Brent-Dekker method which combines an interpolation strategy with the bisection algorithm. Bisection method or interval halving is the simplest bracketing method for root finding of a continuous non-linear function, namely $f(x)$. This method has a linear convergence rate [1]. The first step in bisection method is to provide a search bound. The search bound represented by $[x_a \; x_b]$ is the limit of an interval where the sign of $f(x_a)$ and $f(x_b)$ are different. Based on the Intermediate Value Theorem, when $f(x_a)$ and $f(x_b)$ have opposites signs, then there is at least one real root between $x_a$ and $x_b$. The iteration begins with halving the search space. The midpoint of the interval $x_m = (x_a + x_b)/2$ and $f(x_m)$ are evaluated. If $f(x_a) \times f(x_m) < 0$, then $x_b$ is replaced by $x_m$ otherwise $x_a$ is replaced by $x_m$. The search space is halved at each step. The process is repeated and the root estimate refined by dividing the subintervals into finer increments. Iteration terminates if at any point $f(x)$ equals 0 or the distance between $x_a$ and $x_b$ becomes lower than a requested precision. This well-known algorithm can find only one root in the search bound. If there is more than one root, it is unclear which root is found. It's very favorable to extend bracketing methods in somehow to find all roots of an objective function at a given interval. Such a method is named global root bracketing method (GRBM). In this paper a globalization of bisection method for finding all roots of an objective function is proposed.

---

[1] Vali-e Asr 7, No. 4, Jajarm, Iran, Postal code: 9441639133

The rest of the paper is organized as follows. In section 2, the GRBM which previously proposed will be reviewed and its deficiencies will be shown. In section 3, a new approach which adds adaptive mesh refinement (AMR) feature to the GRBM is introduced. It will be shown that implementing AMR can greatly increase the speed of GRBM, adding new features, and suppressing the deficiencies. Finally, a conclusion is given in Section 4. In this paper, the root finding is discussed for one-dimension problems but the idea can be applied to higher dimensions and will be discussed in further publications.

## 2   Global root bracketing method

Using bracketing methods for finding all the roots of an objective function has previously reported [2-8]. For a detailed review of root bracketing methods refer to [1]. Such methods sometimes called bisection-exclusion methods. In the GRBM, incremental search approach is used to locate subintervals where the function changes sign. These subintervals are called bracketing intervals. Other subintervals which fail to pass bisection criterion are excluded. A bracketing method is used to locate the root in every one of the bracketing intervals. The generic form of GRBM can be presented in the following algorithm [1]:

1. Input: $f(x) = 0$, search bound limits, halving threshold (*HT*), tolerance for stopping criterion
2. Create an initial mesh
3. Evaluate $f(x)$ at the nodes of the initial mesh
4. Select bracketing intervals
5. Apply a root bracketing method to each bracketing interval
6. Output roots founded in step 5

In some works, the search interval divided into sub-intervals all in once [5]. In some other works, dividing goes on step by step [2]. If an interval is not detected as a bracketing interval and its size is bigger than the *HT*, it should be halved. That's because failing in fulfilling bracketing criterion doesn't guarantee that a subinterval does not contain any root. So, it should be chuck down further. The *HT* is chosen to limit halving process. The problem which arises is that a small *HT* is computationally costly. On the other hand, if a large *HT* is chosen, some roots may be omitted. The problem is compounded by the possible existence of even-multiplicity roots, such as $f(x) = (x - x_0)^2$. A representation of GRBM is depicted in Fig. 1.

Dashed lines in Fig. 1 show increment length or mesh. One can see that the *HT* is not small enough to detect bracketing intervals which enclose any of the two roots which lie at [$x_3$ $x_4$]. The last root on the right is a root with even-multiplicity and would be missed regardless of *HT* because the curve never crosses *x* axis at an even-multiple root. Roots lie at [$x_1$ $x_2$] and [$x_2$ $x_3$] can be detected correctly. It's clear from Fig. 1 that GRBM leads to static mesh refinement with constant grid spacing.

Here, the first deficiency of the generic form of GBM is discussed. There are regions in the search space which need to be investigated with smaller *HT* ([$x_4$ $x_5$] in Fig. 1). On the other hand, there are regions ([$x_1$ $x_2$] in Fig. 1) that a large *HT* can take care of root finding. Uniform mesh refinement which is depicted in Fig. 1 fails to find bracketing intervals efficiently. In order to separate close roots with static mesh refinement, *HT* should not be bigger than the minimum distance between roots. In section 2, an adaptive mesh refinement will be presented which can solve this problem.

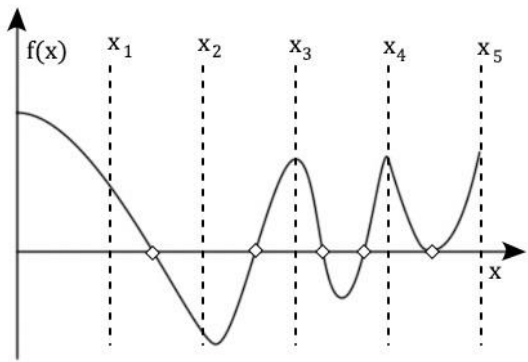

Figure 1: Cases where roots could be missed in GRBM

The other problem is related to definition of error or termination criterion in root bracketing methods. The acceptable form of error for root bracketing methods including the bisection method is same as other iterative methods and can be expressed as follows [9]:

$$\left|\frac{x_m^{new} - x_m^{old}}{x_m^{new}}\right| = \varepsilon_a \qquad (1)$$

in which, $x_m^{new}$ is the root for the present iteration and $x_m^{old}$ is the root from the previous iteration. $\varepsilon_a$ is the tolerance of $x_m^{new}$. The absolute value is used because the sign of the tolerance is not important. When $\varepsilon_a$ becomes less than a pre-specified stopping criterion $\varepsilon_s$, the computation is terminated. To avoid dividing by very small numbers, relative tolerance is presented as follows [9]:

$$\left|\frac{x_m^{new} - x_m^{old}}{\max(1, x_m^{new})}\right| = \varepsilon_a \qquad (2)$$

When there are roots very close to each other, the $\varepsilon_s$ should be lower than the distance between two roots in order to separate roots precisely. But, there is no foreknowledge about how close to each other two distinct roots can be. In fact, that's what we are going to find out. So, decreasing $\varepsilon_s$ cannot be a solution to precise reporting of close roots. Also, very small $\varepsilon_s$ will increase computation time. An adaptive $\varepsilon_s$ should be defined to be decreased in relation to the distance between the roots. An appropriate definition will be presented in the next section which enables GRBM to detect close roots precisely.

## 3     Global root bracketing method with adaptive mesh refinement

The flowchart of GRBM with AMR is presented in Fig. 2. Bisection method is used as the root bracketing method. Each block is numbered to be addressed easily. Inputs which are defined in box 1 include $f(x) = 0$ as the objective function and $x_L(1)$ and $x_R(1)$ as the left and right limits of the search bound. Also, $\varepsilon$, $\varepsilon_m$, and $\varepsilon_f$ which are related to tolerance estimation, and $C$ which is related to $HT$ should be provided by user. In box 2, some counting parameters are initialized. In box 3, the midpoint of subinterval is calculated and stored in a row matrix $x_m(j)$. j is an index number for subintervals. In box 4, the bisection criterion is checked. If it is found out that the subinterval is not a bracketing one, algorithm proceeds to box 5. If the objective function returns a value below $\varepsilon_f$ at $x_m(j)$ (box 5), then the midpoint is saved as a root at box 6. $\varepsilon_f$ should be a number very close to zero. If the condition in box 5 does not meet, the condition in box 7 should be checked. Box 7 states that the subinterval should be halved if its size is bigger than $HT$. Instead of assigning a constant value to $HT$, a new definition for $HT$ is presented as follows:

$$HT = C \times \frac{min\{|f(x_L(j))|, |f(x_R(j))|\}}{x_R(j) - x_L(j)} \qquad (3)$$

$C$ is defined by user. It can adjust how deep the halving process can go. Other parameters in Eq. 3 change as the algorithm proceeds. $C$ is multiplied by the minimum of absolute value of $f(x)$ at each of the left and right limit of the subinterval. Multiplying by such term means that when $f(x)$ is close to zero at any limit of the subinterval, the $HT$ is decreased. That's because it's much more probable to find a root in the vicinity of a point where $f(x)$ is close to zero. On the other hand, when $f(x)$ is far from zero, this chance is decreased. The problem with multiplying with the mentioned term arises when the algorithm converges to a root. In this case, $f(x)$ at both ends of the subinterval becomes so close to zero. This leads $HT$ to become very small step by step. To solve this problem, we should consider that as the bisection converges to a root, the subinterval becomes small. So, by dividing the size of the working interval to $C$, nonstop decreasing of $HT$ is inhibited. In contrast to generic form of GRBM, the definition of $HT$ which was presented in Eq. 3 will lead to an adaptive mesh refinement. In numerical analysis, adaptive mesh refinement is a method of changing the accuracy of a solution in certain regions, during the time the solution is being calculated. In regions where $f(x)$ is close to zero, there would be a dense mesh refinement and in regions where $f(x)$ is far from zero there would be a coarse mesh refinement. The proposed approach decreases function evaluations and computation time considerably.

If the condition in box 7 (Fig. 2) is met, the working interval is halved. The location of left and right limit of new subintervals are saved in row matrices of $x_L$ and $x_R$ (box 8). $i$ index is used to number new subintervals which produced by halving process all through the flowchart. In box 20, $j$ will be increased by one and if the condition in box 21 doesn't met, the process goes back to box 3. In fact, the subinterval is changed and bracketing condition is checked again.

In box 4 (Fig. 2), if the subinterval is detected as a bracketing one, the size of bracketing interval is stored in a parameter called $L$ (box 9). That's because $x_L(j)$ and $x_R(j)$ are going to change in the bisection loop. After that, the tolerance check is performed in box 10. Instead of using the relative form of tolerance presented in Eq. 2, new definition for error is presented as follows:

$$Tol = min\{\varepsilon \times L, \varepsilon_m\} \qquad (4)$$

in which, $Tol$ stands for tolerance and $L$ is the size of the bracketing interval. $Tol$ has a maximum of $\varepsilon_m$ but its minimum is boundless. Actually, the minimum of $Tol$ is dependent on the $L$. This type of tolerance definition can solve the problem of reporting close roots with appropriate precision. When two roots are very close to each other, the bracketing interval for detecting each one of them becomes very small in comparison to far a parted roots. In this cases, by multiplying the size of bracketing interval ($L$) into $\varepsilon$, tolerance will be decreased in proportion to the closeness of the roots. In cases where L becomes bigger than one, the minimum operator will choose $\varepsilon_m$ as the tolerance.

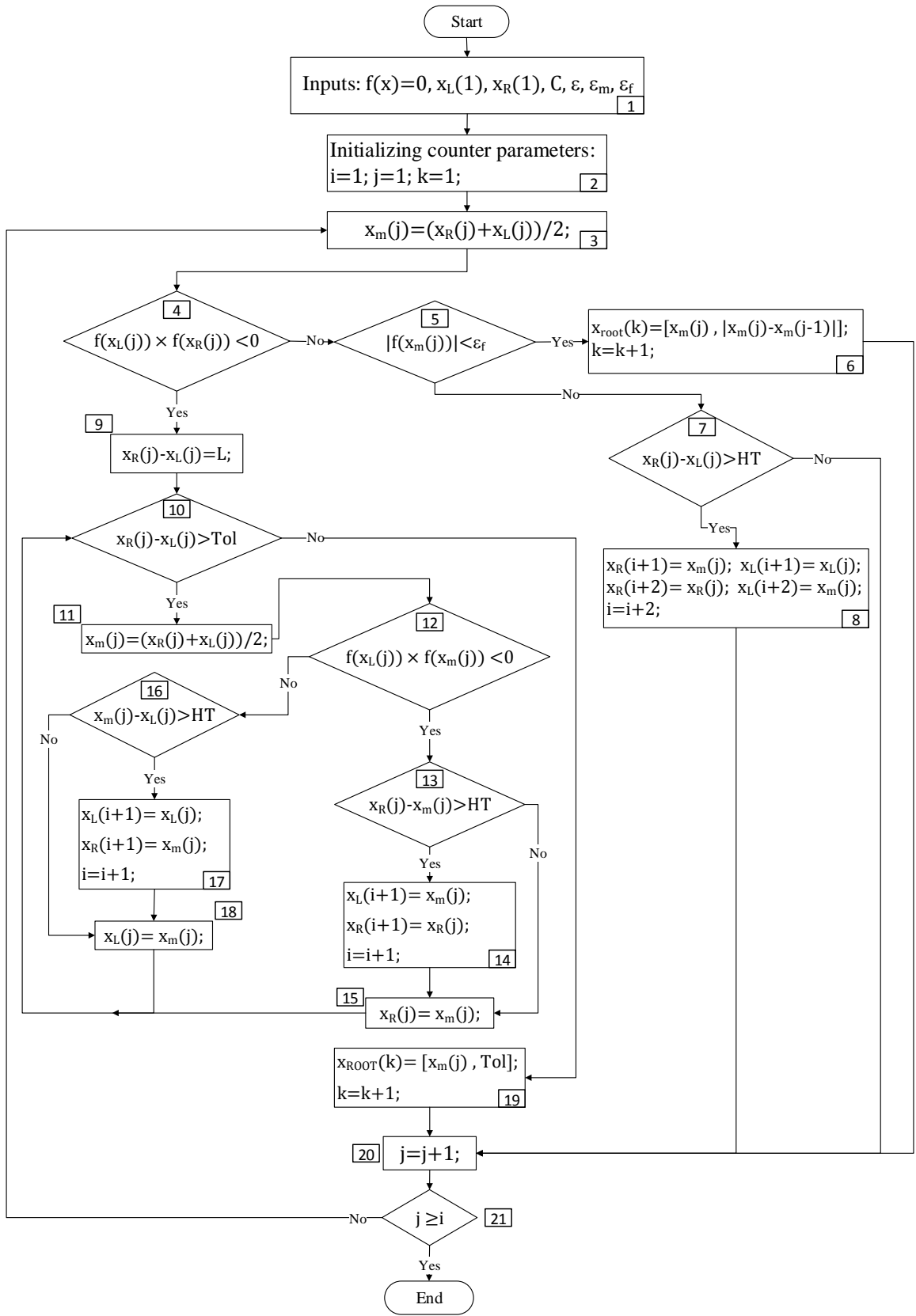

Figure 2: Flowchart of global root bracketing method with adaptive mesh refinement

If the condition in box 10 (Fig. 2) is met, the working interval is halved. One of the halves is chosen as bracketing interval and the other is checked weather if it's bigger than $HT$ or not. If it is found that the non-bracketing half is bigger than $HT$, its location will be saved as a new subinterval which have to be fed into the algorithm afterward. These steps take place in boxes 11 to 18. New subintervals are defined in boxes 8, 14, and 17. Bisection loop is terminated when the condition in box 10 is failed to meet. After failing at box 10, the location of midpoint of the bracketing interval is saved as a root in box 19. In box 20, $j$ is increased by one. If $j$ is lower than $i$, it means that there are still subintervals that should be checked. Location of roots and the error in locating each one them are saved in a matrix called $x_{ROOT}$ (boxes 5 and 19).

It should be noted that the tolerances used in the flowchart cannot be less than machine epsilon of the processor in which the algorithm is executed. Machine epsilon gives an upper bound on the relative error due to rounding in floating point arithmetic. So, roots which are close to each other less than machine epsilon cannot be detected. For most of computing software, machine epsilon is equal to $2.22 \times 10^{-16}$.

## 3.1 Evaluating the GRBM with AMR

In this section the power of the proposed GRBM with AMR is evaluated and compared with the generic form of GRBM through some examples. A program based on flowchart presented in Fig. 2 was written in Matlab software. The objective function was written in a different file within the Matlab environment. In order to minimize function evaluations, the function is defined in such a way that inhibit recalculation. The number of function evaluations is counted within the function file. As an example, a polynomial is chosen for root finding which is as follows:

$$f(x) = (x - 0.5)(x - 0.50001)(x - 4)(x - 4.05)(x - 9.3) \quad 0<x<10 \tag{5}$$

It's clear that $f(x)$ in Eq. 5 has five roots including $x_1$=0.5, $x_2$=0.50001, $x_3$=4, $x_4$=4.05, $x_5$=9.3. $x_1$ and $x_2$ are very close, $x_3$ and $x_4$ are relatively close, and $x_5$ is far apart from other roots. To quantify the concept of closeness, an index of closeness is defined as follows:

$$CI = \frac{the\ size\ of\ initial\ search\ bound}{the\ minimum\ of\ the\ distance\ between\ two\ roots} \tag{6}$$

$CI$ is the closeness index which shows how close to each other two distinct roots can be in relation to the size of the initial search bound. A root with higher $CI$ is harder to find in comparison to a root with lower $CI$. Based on the definition presented in Eq. 6, roots of the polynomial in Eq. 5 has the following $CI$, $CI_1$=$CI_2$=$10^6$, $CI_3$=$CI_4$=100, and $CI_5$=2.04. So, roots cover a wide range of $CI$ which makes the polynomial of Eq. 5 a challenging test for a GRBM.

One of the most important parameters in GRBM with AMR is $C$ in Eq. 3. Other parameters which should be set are $\varepsilon$, $\varepsilon_m$, and $\varepsilon_f$. Program which was wrote in Matlab ran with $\varepsilon_f = 2.22 \times 10^{-16}$ and different values of $C$, $\varepsilon$, $\varepsilon_m$. With every different value for $C$, different number of function evaluations is needed for solving the problem. The distance between points where the objective function is evaluated reveals how the mesh refinement is performed. Table 1 shows relevant data of each run.

Every raw in table 1 is related to a different run which are numbered consecutively in column one. $N$ in column five from the left is the number of function evaluations for each run. In run #1, not all the roots are found. In fact, two roots with highest *CI* are missing. Tolerance for $x_3$ and $x_4$ is almost the half of tolerance for $x_5$. That's because $x_3$ and $x_4$ are close to each other. In run #2, $C$ is lowered to 0.04 which lead to finding all the roots. The minimum tolerance for separating $x_1$ and $x_2$ is $10^{-5}$. In this run, the error for $x_1$ and $x_2$ is $10^{-7}$ which is good enough. To give a picture how differing other parameters affect N, other runs are listed in table 1.

Table 1: Roots of $f(x)$ defined in Eq. 5

| No. | $C$ | $\varepsilon$ | $\varepsilon_m$ | N | Roots |
|---|---|---|---|---|---|
| 1 | 0.05 | $10^{-2}$ | $10^{-3}$ | 45 | $9.300\pm10^{-3}$, $4.0000\pm4\times10^{-4}$, $4.0501\pm4\times10^{-4}$ |
| 2 | 0.04 | $10^{-2}$ | $10^{-3}$ | 95 | $9.300\pm10^{-3}$, $4.0000\pm4\times10^{-4}$, $4.0501\pm4\times10^{-4}$, $0.5000000\pm10^{-7}$, $0.5000099\pm10^{-7}$ |
| 3 | 0.04 | $10^{-2}$ | $10^{-5}$ | 111 | $9.29999\pm10^{-5}$, $4.00000\pm10^{-5}$, $4.05000\pm10^{-5}$, $0.5000000\pm10^{-7}$, $0.5000099\pm10^{-7}$ |
| 4 | 0.01 | $10^{-2}$ | $10^{-5}$ | 139 | $9.29999\pm10^{-5}$, $4.00000\pm10^{-5}$, $4.05000\pm10^{-5}$, $0.5000000\pm10^{-7}$, $0.5000099\pm10^{-7}$ |
| 5 | 0.01 | $10^{-4}$ | $10^{-5}$ | 157 | $9.29999\pm10^{-5}$, $4.00000\pm4\times10^{-6}$, $4.05000\pm4\times10^{-6}$, $0.5000000000\pm10^{-10}$, $0.5000010000\pm10^{-10}$ |

Plots for different values of $C$ which show *HT* and $f(x)$ versus $x$ based on function evaluations which were performed in the program are presented in Fig. 3. Fig. 3a and b are related to run #1 and Fig. 3c and d are related to run #2. It can be seen from Fig. 3b and d that with decreasing *HT*, roots with higher *CI* which are $x_1$ and $x_2$ can be detected. Fig. 3a and c show how *HT* changes as $x$ approaches a root or get far from it. It's clear that the GRBM with AMR can optimize *HT* to minimize function evaluations, enhances the GRBM, and decreases computation time.

To find roots of polynomial presented in Eq. 2 with generic form of GRBM which was summarized in section 2, the algorithm presented in Fig. 2 can be used with two replacement. Replacements are the definition of *HT* and *Tol*. For generic form of GRBM, tolerance is defined as it presented in Eq. 2 and *HT* is regarded as a constant number which should be defined by user. The *HT* should not be bigger than the minimum distance between two roots otherwise some roots will be missed. The minimum distance between roots for problem defined in Eq. 5 is between $x_1$ and $x_2$ and it's equal to $10^{-5}$. So, *HT* should set to be $10^{-5}$ at least. At this value of *HT*, $f(x)$ should be calculated in $10^6$ points. The number comes from initial search bound width divided by *HT*. In fact, the number of function evaluations has the order of magnitude of the highest closeness index between roots. $10^6$ times function evaluations finds bracketing intervals which contain all the roots. According to table 1, EGB method only needs 95 times function evaluations to find all the roots. The difference between GRBM with AMR and generic form of GRBM method in terms of speed and efficiency is substantial. This difference becomes even greater as the *CI* increases.

## 3.2 Finding roots with even-multiplicity

It was noted in section 1 that root bracketing methods are unable to find roots with even-multiplicity because the function curve never crosses x axis in such points but only touches it. AMR embedded in the algorithm presented in Fig. 2 enables us to find roots with even-multiplicity. Although the condition in box 6 always fails in such points but even multiple roots are caught in box 4. GRBM with AMR can converge to roots with even-multiplicity because of the definition proposed for $HT$ in Eq. 3. The $HT$ decreases in the vicinity of such points and lead the algorithm to converge in box 4. As an example, a function with even-multiple roots is defined as follows:

$$g(x) = (x-3)^2(x-4)^2 \quad 0 < x < 5 \tag{7}$$

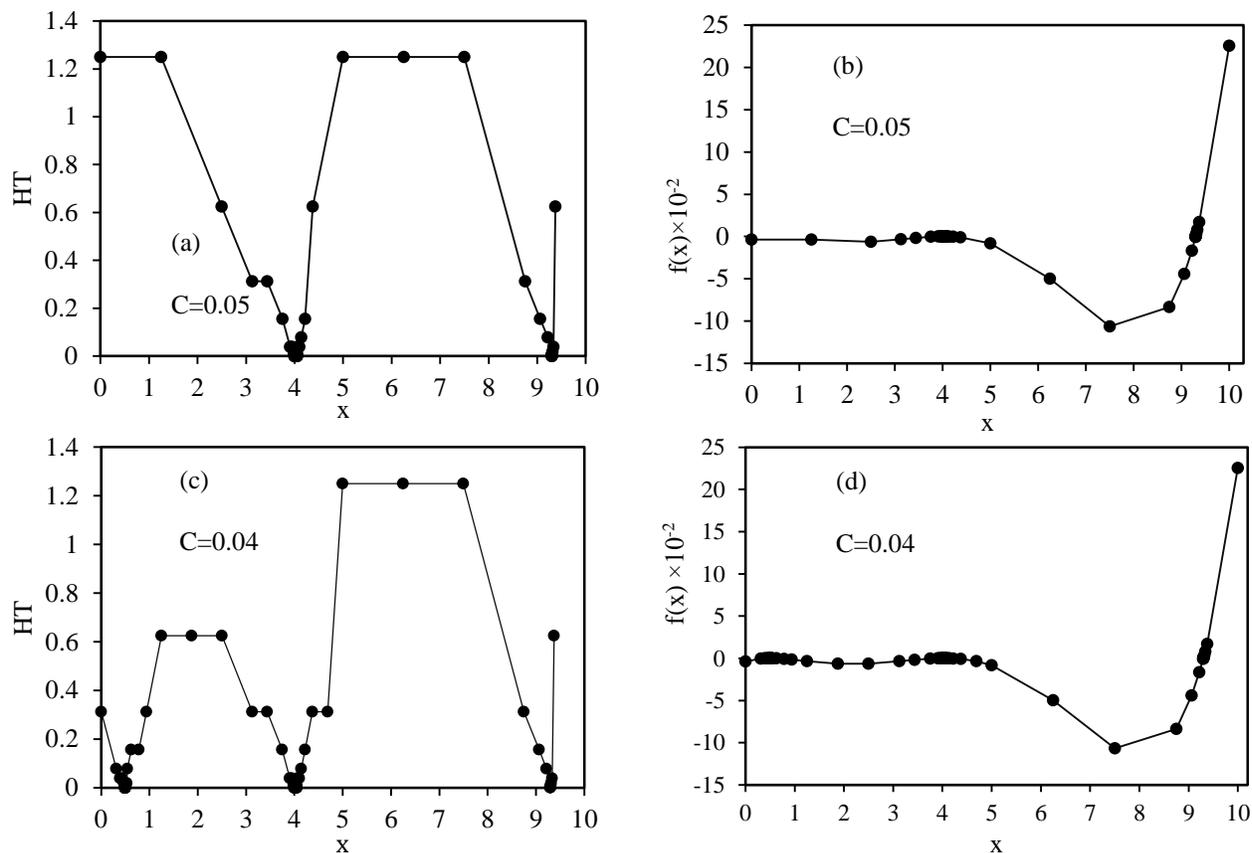

Figure 3: Plots of $HT$ and $f(x)$ versus $x$ for different values of $C$

At $x = 3$ and $4$, $g(x)$ meets its roots which both of them are double roots. Table 2 shows information of different runs for finding roots of $g(x)$ using GRBM with AMR. The value of $\varepsilon$ and $\varepsilon_m$ are irrelevant because the process of root finding never enters bisection loop in this case. But, the parameters $C$ and $\varepsilon_f$ does matter. In all of the runs, $C$=4. The algorithm was executed on a processor in which the floating point format has double precision with machine epsilon of $2.22\times10^{-16}$.

Plot of g(x) in run 1 (table 2) is depicted in Fig. 4. It's clear from Fig. 4 and table 2 that roots of $g(x)$ are detected properly. Because even-multiple roots are detected in a different box (Fig. 2), they can be easily marked as even-multiple.

It's clear from table 2 that the lowest achievable tolerance for roots of g(x) is $10^{-7}$ for $\varepsilon_f$ equal to machine epsilon. What if one wants to achieve a lower tolerance? To do this, a more strict condition need to be applied. If the objective function is differentiable, the first derivative at a point with even multiple root should be equal to zero. This fact is shown in Fig. 1 at the last root on the right. So, to get more precise result from the algorithm proposed in Fig. 2, the condition in box 4 is removed and the following code is implemented:

```
IF |f(x_m(j))| < ε_f
    IF |df(x_m(j))| < ε_d
        Output x_m(j) as an even-multiple root.
    END
END
```

Table 2: Roots of $g(x)$ defined in Eq. 7

| No. | $\varepsilon_f$ | N | Roots |
| --- | --- | --- | --- |
| 1 | $10^{-8}$ | 50 | $3.000\pm10^{-3}$, $3.999\pm10^{-3}$ |
| 2 | $2.22\times10^{-16}$ | 98 | $3.0000000\pm10^{-7}$, $4.0000000\pm10^{-7}$ |

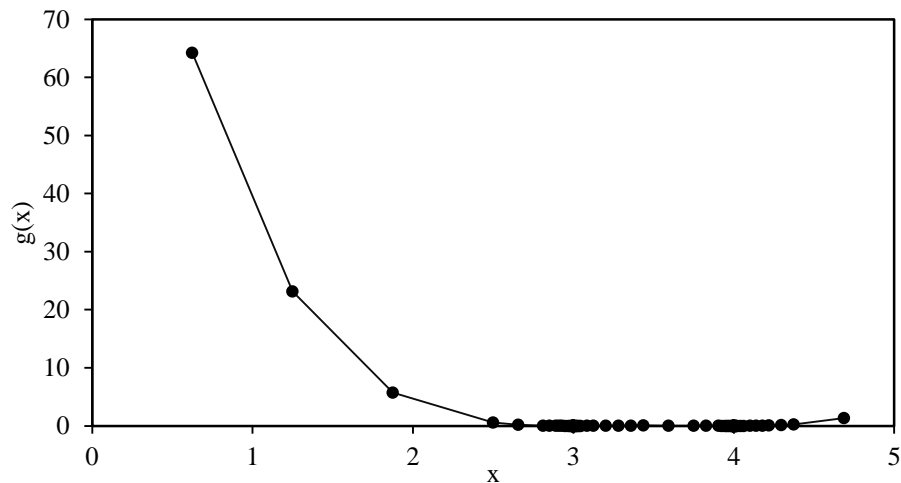

Figure 4: plot of g(x) at run 1 described in table 2.

in which, $df$ is the first derivative of $f$ and $\varepsilon_d$ is the tolerance for zero. In this way, both of the function and its derivative should be checked. For getting the most accurate result, both of $\varepsilon_f$ and $\varepsilon_d$ should be equal to machine epsilon. The program used for finding roots of g(x) ran with the proposed codes. $\varepsilon_f$ and $\varepsilon_d$ set to machine epsilon of $2.22\times10^{-16}$. Result show that roots are found properly with accuracy of $2.22\times10^{-16}$. 202 times function evaluations and 198 times first derivative evaluations were performed to detect roots with the highest possible accuracy.

## 3.3 Finding roots with odd-multiplicity

Generally speaking, root bracketing methods can find roots with odd-multiplicity but at a reduced rate of convergence. Methods with higher rate of convergence like Newton's method are used to find roots in such cases. In this section, finding roots with odd-multiplicity using the GRBM with AMR is discussed.

A sample objective function which contains a triple root is defined as follows:

$$h(x) = (x - 0.5)^3(x - 0.50001)(x - 1) \quad 0 < x < 1.5 \tag{8}$$

$h(x)$ has three roots as follows, $x_1$=0.5, $x_2$=0.50001, and $x_3$=1. GRBM with AMR is used to find roots of $h(x)$. It was found out that the algorithm with $\varepsilon = \varepsilon_m = 10^{-5}$ and $C$=500 needs at least 500 times function evaluations to find the roots. But after finding the three roots, the algorithm continues to divide subintervals into halves around the triple root consistently. This increases the number of function evaluations into millions. In fact, the algorithm get stuck between near the triple root. This phenomena doesn't happen for simple close roots as it was shown in section 3.1. If the root $x_2$, instead of being so close to $x_1$, is moved to the point 0.6; after 1491 times function evaluations, all the roots are found. So, finding odd-multiple roots increase number of function evaluations drastically and when they are so close to another root makes the problem even bigger. The reason behind this drastic increase of function evaluations for finding odd-multiple roots lies in this fact that the objective function crosses the x axis with a flat shape and very smoothly in such points. So, in a vicinity around an odd-multiple root, objective function stays very close to x axis. This behavior deceives the algorithm to search deeper and deeper to find another root where there is none. Adding another root close to an odd-multiple root makes this case harder.

To solve the problem of finding odd-multiple roots, a modification in the formula of *HT* which defined in Eq. 3 is done as follows:

$$HT = C \times \frac{min\{|f(x_L(j))|, |f(x_R(j))|\}}{(x_R(j) - x_L(j))^n} \tag{9}$$

The difference between Eq. 3 and 9 is that the denominator on the right side of formula is raised to power of $n$. The proposed definition for *HT* is used to find roots of Eq. 8. The result shows that with $\varepsilon = \varepsilon_m = 10^{-5}$, $C$=20, and $n$=3 only 87 times function evaluations is enough to detect roots as follows, $x_1$=0.4999999999±10⁻¹⁰, $x_2$=0.5000099999±10⁻¹⁰, $x_3$=1.00000±10⁻⁵. The right value of $n$ and $C$ for minimum number of function evaluations are obtained by trial and error. For subintervals smaller than 1 in width, as the n increases for $n$>1, the *HT* increases. This is how the excess halving process is inhibited. Adjusting the right value of $n$ can greatly decrease function evaluations. But, there is a serious drawback with increasing $n$ higher than 1. In cases where simple roots especially close ones are distributed in the initial search bound, $n$>1 lead to missing some of simple roots.

## 3.4 Finding roots in the most difficult cases

To have a perfect algorithm, $n$ in Eq. 9 have to change based on how the function behaves at different regions of the search bound. In the absence of such a definition for $n$, an alternative strategy for finding all the roots when there are combination of different types of root is proposed here. The strategy relies on running the algorithm different times with different values of $n$. In this way, at $n$>1 multiple roots are found with the least effort. At the next run, regions which include the roots are excluded from the search domain and the rest of search interval fed into algorithm with $n=1$. The limits of these regions should be defined a little far from roots to avoid the algorithm to converge toward the excluded roots in the next run. This strategy can find all the roots regardless of the type and closeness index of roots in an efficient way. As an example, roots of an objective function which is formulated as followed are found by the mentioned strategy.

$$p(x) = (x - 0.5)^3(x - 0.50001)^3(x - 4.0)(x - 4.0001)(x - 4.2)^2 \qquad 0 < x < 4.5 \qquad (10)$$

$p(x)$ has five roots as follows, $x_1=0.5$, $x_2=0.50001$, $x_3=0.4$, $x_4=0.401$, and $x_5=4.2$. It's clear that different types of roots are placed close together which makes the problem of root finding very hard. In the first run, $n=5$, $C=0.1$, and $\varepsilon = \varepsilon_m = 10^{-5}$. Also, in the first run, because we look only for odd-multiple roots, box 5 in Fig. 2 which is for detecting even-multiple roots is removed. After 78 times function evaluations, roots are found as follows, $x_1=0.4999999999\pm10^{-11}$ and $x_2= 0.5000099999\pm10^{-11}$. At the second run, $n=1$, $C=0.01$, $\varepsilon = \varepsilon_m = 10^{-5}$, and the condition on box 5 in Fig. 2 comes back with the adjustment that proposed in section 3.2. In the second run, the search domain is replaced by [0.6 4.5]. After 292 times function evaluations and 125 times first derivative evaluations, the rest of the roots are found as follows, $x_3= 3.999999999\pm10^{-10}$, $x_4= 4.000100000\pm10^{-10}$, and $x_5= 4.200000000000000\pm10^{-16}$. Use of first derivative information is not compulsory, but such information is used to decrease the error bound for even-multiple roots. Without use of first derivative information, the number of function evaluations decreases to 227 and the error bound for the even-multiple root increases to $10^{-9}$.

Although, in examples which were discussed, polynomials were chosen as the objective function; but methods which proposed in this paper are not dependent on the type of the objective function.

## 4   Conclusion

Some of the most troublesome deficiencies of bracketing methods for root finding were addressed and resolved. New definition of halving threshold was presented which lead to adaptive mesh refinement. The proposed method can decrease function evaluations notably in comparison to global root bracketing with static mesh refinement. In a case study revealed that for a polynomial root finding problem, 95 times function evaluations in global root bracketing with adaptive mesh refinement is corresponding to at least $10^6$ times function evaluations with static mesh refinement. Another achievement is presenting a new definition for stopping criterion of root bracketing methods. This definition solves the problem of reporting close roots with appropriate accuracy. Special cases of finding odd or even-multiple roots were discussed. With adjustment in the definition of halving threshold, the algorithm was enhanced to handle finding multiple roots with the least number of function evaluations. The methods proposed in this paper showed that they can be very successful in global root finding especially when it comes to separating close roots.

References

[1] J.M. McNamee, V.Y. Pan, Chapter 7 - Bisection and Interpolation Methods, in: J.M. McNamee, V.Y. Pan (Eds.) Numerical Methods for Roots of Polynomials - Part II, Elsevier, 2013, pp. 1.
[2] X. Ying, I.N. Katz, A simple reliable solver for all the roots of a nonlinear function in a given domain, Computing, 41 (1989) 317-333.

[3] A. Eiger, K. Sikorski, F. Stenger, A Bisection Method for Systems of Nonlinear Equations, ACM Trans. Math. Softw., 10 (1984) 367-377.
[4] R.B. Kearfott, Some Tests of Generalized Bisection, ACM Trans. Math. Softw., 13 (1987) 197-220.
[5] D. Bachrathy, G. Stépán, Bisection method in higher dimensions and the efficiency number, Period. Polytech. Mech., 56 (2012) 81.
[6] H.S. Wilf, A Global Bisection Algorithm for Computing the Zeros of Polynomials in the Complex Plane, J. ACM, 25 (1978) 415-420.
[7] G.R. Wood, Multidimensional bisection applied to global optimisation, Comput. Math. Appl., 21 (1991) 161-172.
[8] D.J. Kavvadias, F.S. Makri, M.N. Vrahatis, Efficiently Computing Many Roots of a Function, Siam J. Sci. Comput., 27 (2005) 93.
[9] E.K.P. Chong, S.H. Zak, An introduction to optimization, 4rd ed., John Wiley & Sons, Hoboken, New Jersey, 2013.